\documentclass[10 pt]{amsart}

\usepackage[colorlinks=true,linkcolor=blue,citecolor=blue,urlcolor=blue,hypertexnames=false]{hyperref}
\usepackage{bookmark}
\usepackage{amsthm,thmtools,amssymb,amsmath,amscd}

\usepackage{fancyhdr}
\usepackage{esint}
\usepackage{enumerate}

\usepackage{pictexwd,dcpic}

\swapnumbers
\declaretheorem[name=Theorem,numberwithin=section]{thm}
\declaretheorem[name=Remark,style=remark,sibling=thm]{rem}
\declaretheorem[name=Lemma,sibling=thm]{lemma}
\declaretheorem[name=Proposition,sibling=thm]{prop}
\declaretheorem[name=Definition,style=definition,sibling=thm]{defn}
\declaretheorem[name=Corollary,sibling=thm]{cor}

\declaretheorem[name=Example,style=remark,sibling=thm]{example}

\numberwithin{equation}{section}

\usepackage{cleveref}
\crefname{lemma}{Lemma}{Lemmata}
\crefname{prop}{Proposition}{Propositions}
\crefname{thm}{Theorem}{Theorems}
\crefname{cor}{Corollary}{Corollaries}
\crefname{defn}{Definition}{Definitions}
\crefname{example}{Example}{Examples}
\crefname{rem}{Remark}{Remarks}
\crefname{assum}{Assumption}{Assumptions}
\crefname{notation}{Notation}{Notation}

\usepackage{autonum}

\newcommand{\ti}{\tilde}

\newcommand{\cn}{\colon}
\newcommand{\sub}{\subset}

\newcommand{\N}{\mathbb{N}}

\newcommand{\R}{\mathbb{R}}

\newcommand{\8}{\infty}

\newcommand{\al}{\alpha}

\newcommand{\ka}{\kappa}
\newcommand{\la}{\lambda}

\newcommand{\s}{\sigma}

\newcommand{\Om}{\Omega}

\newcommand{\G}{\Gamma}

\newcommand{\La}{\Lambda}

\newcommand{\cL}{\mathcal{L}}

\newcommand{\cW}{\mathcal{W}}

\newcommand{\cD}{\mathcal{D}}
\newcommand{\cB}{\mathcal{B}}

\newcommand{\del}{\partial}

\newcommand{\fa}{\forall}

\newcommand{\fr}[2]{\frac{#1}{#2}}
\newcommand{\x}{\times}

\DeclareMathOperator{\id}{id}

\DeclareMathOperator{\tr}{tr}

\DeclareMathOperator{\ad}{ad}

\newcommand{\pf}[1]{\begin{proof} #1 \end{proof}}
\newcommand{\eq}[1]{\begin{equation}\begin{alignedat}{2} #1 \end{alignedat}\end{equation}}

\newcommand{\br}[1]{\left(#1\right)}

\newcommand{\ra}{\rightarrow}

\newcommand{\mt}{\mapsto}

\newcommand{\mc}{\mathcal}

\newcommand{\mrm}{\mathrm}

\newcommand{\hp}{\hphantom}
\newcommand{\q}{\quad}

\begin{document}

\title{Isotropic functions revisited}
\author{Julian Scheuer}
\address{Albert-Ludwigs-Universit\"{a}t,
Mathematisches Institut, Eckerstr. 1, 79104
Freiburg, Germany}
\email{julian.scheuer@math.uni-freiburg.de}
\date{\today}


\begin{abstract}
To a real $n$-dimensional vector space $V$ and a smooth, symmetric function $f$ defined on the $n$-dimensional Euclidean space we assign an associated operator function $F$ defined on linear transformations of $V$. $F$ shall have the property that, for each inner product $g$ on $V$, its restriction $F_{g}$ to the subspace of $g$-selfadjoint operators is the isotropic function associated to $f$. This means that it acts on these operators via $f$ acting on their eigenvalues. We generalize some well known relations between the derivatives of $f$ and each $F_{g}$ to relations between $f$ and $F$, while also providing new elementary proofs of the known results. By means of an example we show that well known regularity properties of $F_{g}$ do not carry over to $F$.
\end{abstract}

\maketitle

\section{Introduction}
Consider a function $f\in C^{\8}(\R^{n})$ which is {\it{symmetric}}, i.e.
\eq{f(\ka_{1},\dots,\ka_{n})=f(\ka_{\pi(1)},\dots,\ka_{\pi(n)})\q\fa \pi\in\mc{P}_{n},}
where $\mc{P}_{n}$ is the permutation group on $n$ elements. Let $V$ be a real, $n$-dimensional vector space and $\cL(V)$ be the vector space of linear operators on $V$. If $V$ carries an inner product $g$, on the vector subspace $\Sigma_{g}(V)\sub \cL(V)$ of $g$-selfadjoint operators one can define a map 
\eq{F_{g}\cn \Sigma_{g}(V)&\ra \R\\
			A&\mt f(\mrm{EV}(A)),}
where $\mrm{EV}(A)=(\ka_{1},\dots,\ka_{n})$ denotes the ordered $n$-tuple of real eigenvalues of $A$. In \cite{Ball:/1984} J. Ball proved that if $f\in C^{r}(\R^{n})$, $r=1,2,\8$, the function $F_{g}$ is also of class $C^{r}$. Furthermore, using Schauder theory, he showed that if $f\in C^{r,\al}(\R^{n})$, $r\in \N$, $0<\al<1$, then also $F$ is in the respective function class. Also compare \cite[Sec.~2.1]{Gerhardt:/2006} for a detailed proof of these regularity results.  
For $r\geq 3$, the implication
\eq{f\in C^{r}(\R^{n})\q\Rightarrow \q F_{g}\in C^{r}(\Sigma_{g}(V))}
was proven in \cite{Silhavy:/2000}.

In these results one always starts with an inner product space $(V,g)$. In many applications one has to deal with a whole family of such spaces, where $g$ may vary. For example in geometric curvature problems one is often faced with a map $F$ being evaluated on the Weingarten tensor $\cW$, an endomorphism field with values in the tensor bundle of linear transformations of the tangent spaces. From point to point, these linear maps $\cW(x)$ are self-adjoint with respect to different metrics, so one has to be careful with the domain of $F$.

One may observe, that for the most natural symmetric functions, e.g. 
\eq{s_{1}=\sum_{i=1}^{n}\ka_{i}\q\text{or}\q s_{n}=\prod_{i=1}^{n}\ka_{i}}
there is no ambiguity about how to define $F$ even on the whole space $\cL(V)$ and not only on some $\Sigma_{g}(V)$. Namely for $s_{1}$ just set
\eq{F(A)=S_{1}(A)=\tr(A)}
and for $s_{n}$ set
\eq{F(A)=S_{n}(A)=\det(A).}
The functions $s_{1}$ and $s_{n}$ are special cases of the {\it{elementary symmetric polynomials}} $s_{k}$, $1\leq k\leq n$, cf. \cref{skpk}, to which we associate
\eq{S_{k}(A)=\fr{1}{k!}\fr{d^{k}}{dt^{k}}\det(I+tA)_{|t=0}.}
 It is true that every symmetric function $f\in C^{\8}(\G)$ on a symmetric open set $\G\sub\R^{n}$ can be written as a function of the $s_{i}$,
\eq{f=\rho(s_{1},\dots,s_{n}),}
where $\rho\in C^{\8}(U)$ for some open $U\sub \R^{n}$, cf. \cite{Glaeser:01/1963}. In case $f\in C^{r}(\G)$, $\rho$ will in general have less regularity, cf. \cite{Barbancon:/1972}. In both cases the function 
\eq{F=\rho(S_{1},\dots,S_{n})}
is defined on an open set $\Om\sub\cL(V)$ and satisfies
\eq{F(A)=f(\mrm{EV}(A))}
for all $\R$-diagonalisable $A\in \cL(V)$ with eigenvalues in $\G$. Hence $F$ can be differentiated in all directions of $\cL(V)$.

The aim of this short note is a transfer of some well known and often used relations between derivatives of $F$ and $f$ to the new situation, that $F$ can be differentiated in all of $\cL(V)$. In previous treatments of this, only the relation between $f$ and $F_{g}$ was studied for some fixed metric $g$, compare for example \cite{Andrews:/2007,Ball:/1984,BowenWang:01/1970,ChadwickOgden:01/1971,ChadwickOgden:01/1971b,Gerhardt:/2006,HuiskenSinestrari:09/1999,LewisSendov:/2001,Sendov:07/2007,Silhavy:/2000}. Our approach is by direct calculation of the proposed relations for the elementary symmetric polynomials and then to transfer them to general functions. Note that this approach also provides a new, quite elementary proof of the corresponding results for the pair $(f,F_{g})$ with fixed inner product $g$, as obtained in \cite[Thm.~5.1]{Andrews:/2007} and \cite[Lemma~2.1.14]{Gerhardt:/2006}. 

The motivation to write this note came up during the preparation of \cite{BIS4}, where we had to apply derivatives of $F_{g}$ to some non-$g$-selfadjoint operators, so the need for a globally defined $F$ was apparent. For illustration, have a look at the following simple example:

\begin{example}\label{example}Let $f$ be the second power sum,
\eq{f(\ka)=\sum_{i=1}^{n}\ka_{i}^{2},\quad F(A)=\tr(A^{2}),}
then $F$ is clearly the associated operator function for $f$ and $F$ is defined on whole $\cL(V)$. $f$ is a convex function of the $\ka_{i}$. However, 
\eq{F\cn \cL(V)\ra \R}
is {\it{not}} convex: Indeed there holds
\eq{dF(A)B=2\tr(A\circ B),}
\eq{d^{2}F(A)(B,C)=2\tr(B\circ C)}
and hence
\eq{d^{2}F(A)(\eta,\eta)=2\tr(\eta^{2})<0}
for a nonzero skew-symmetric (with respect to a basis of eigenvectors of $A$) $\eta$.
\end{example}

 The fact that $F$ is in general not convex, when considered as a function on $\cL(V)$, caused trouble in the preparation of \cite{BIS4}, where we had to estimate the term $d^{2}F(\dot{\cW},\dot{\cW})$ along some curvature flow. Here $\dot{\cW}$ is the evolution of the Weingarten tensor. For the particular flow considered in \cite{BIS4}, we could not prove the symmetry of $\dot{\cW}$. This trouble was the main motivation to write this note and to extend the formulas for derivatives of $F$, as in \cref{2DerHk}.

\section{Symmetric functions and associated operator functions}\label{Poly}

For an $n$-dimensional, real vector space $V$, the aim of this section is to deduce relations between the derivatives of the functions $f$ and $F$ as described in the introduction. First we fix some definitions and notation.

{\defn{\label{skpk}
On $\R^{n}$ we denote the elementary symmetric polynomials for $1\leq k\leq n$ by $s_{k},$
\eq{s_{k}(\ka_{1},\dots,\ka_{n}):=\sum_{1\leq i_{1}<\dots<i_{k}\leq n}\prod_{j=1}^{k}\ka_{i_{j}}}
and the power sums for all $ k\in \N$ by $p_{k},$
\eq{p_{k}(\ka)=\sum_{i=1}^{n}\ka_{i}^{k}.}
}}
{\defn{(i)~Let $V$ be an $n$-dimensional real vector space and $\cD(V)\sub \cL(V)$ be the set of diagonalisable endomorphisms. Then we denote by $\mrm{EV}$ the eigenvalue map, i.e.
\eq{\mrm{EV}\cn \cD(V)&\ra \R^{n}/\mc{P}_{n}\\
				A&\mt (\ka_{1},\dots,\ka_{n}),}
where $\ka_{1},\dots, \ka_{n}$ denote the eigenvalues of $A$ and $\mc{P}_{n}$ is the permutation group of $n$ elements.

(ii)~Let $\G\sub\R^{n}$ be open and symmetric, then we define
\eq{\cD_{\G}(V)=\mrm{EV}^{-1}(\G/\mc{P}_{n}).} 
}}

{\rem{Note that $\mrm{EV}$ is continuous, compare \cite{Zedek:/1965}. }}

{\lemma{\label{Hk} Let $V$ be an $n$-dimensional real vector space. Then for all $k\in \N$ there exists a function $P_{k}\in C^{\8}(\cL(V))$ with
\eq{\label{Pk-A}P_{k}(A)=p_{k}\circ \mrm{EV}(A)\q\fa A\in \cD(V).}  }}

\pf{
 Simply set
\eq{P_{k}(A)=\tr(A^{k}).}
Then there holds
\eq{P_{k}(A)=p_{k}(\mrm{EV}(A))\q\fa A\in \cD(V).}
}

Since the $P_{k}$ are smooth, we want to investigate the structure of their derivatives.

{\prop{\label{DerHk} Let $V$ be an $n$-dimensional real vector space. Let $U\sub\R^{m}$ be open and $\psi\in C^{r}(U)$, $r\geq 1$. Then the function
\eq{f=\psi(p_{1},\dots,p_{m})}
is defined on an open symmetric set $\G\sub\R^{n}$ and the function $F=\psi(P_{1},\dots,P_{m})$ is defined on an open set $\Om\sub\cL(V)$. There holds  \eq{F_{|\cD_{\G}(V)}=f\circ\mrm{EV}_{|\cD_{\G}(V)}}
 and the derivatives of $F$ evaluated at a fixed $A\in \Om$ are given by
\eq{\label{DerHk-4}dF(A)B=\tr(F'(A)\circ B)=\sum_{l=1}^{m}l\fr{\del\psi}{\del P_{l}}\tr\br{A^{l-1}\circ B}\q\fa B\in \cL(V),}
where
\eq{\label{DerHk-1}F'(A)=\sum_{l=1}^{m}l\fr{\del \psi}{\del P_{l}}A^{l-1}.} }}

\pf{
Only the formula for $dF$ has to be checked, while all other statements are obvious. The function $P_{1}(A)=\tr(A)$ is linear and hence
\eq{dP_{1}(A)B=\tr(B)\q\fa A,B\in \cL(V).}
Furthermore by the chain rule there holds
\eq{\label{DerHk-3}dP_{k}(A)B=d(P_{1}(A^{k}))(A)B=k\tr(A^{k-1}\circ B)\q\fa A,B\in \cL(V).}
Thus
\eq{\label{DerHk-2}dF(A)B=\sum_{l=1}^{m}\fr{\del \psi}{\del P_{l}}dP_{l}(A)B=\tr(F'(A)\circ B)}
and hence the proof is complete. }

{\rem{It is well known that the elementary symmetric polynomials $s_{k}$ are functions of the $p_{k}$, cf. \cite{Mead:10/1992}, and hence \cref{DerHk} also applies to these. }}

{\cor{\label{EVDerHk}Let $V$ be an $n$-dimensional real vector space and let $f$ and $F$ be as in \cref{DerHk}. Suppose $A\in \cD_{\G}(V)$. Then the endomorphisms $F'(A)$ and $A$ are simultaneously diagonalisable. For a basis $(e_1,\dots,e_n)$ of eigenvectors for $A$ with eigenvalues $\kappa=(\ka_1,\dots,\ka_n),$ the eigenvalue $F^i$ of $F'(A)$ with eigenvector $e_i$ is given by
\eq{\label{EVDerHk-1}F^{i}(A)=\fr{\del f}{\del\ka_{i}}(\ka).} 
}}

\pf{That $F'(A)$ and $A$ are simultaneously diagonalisable follows from \eqref{DerHk-1} immediately. Let $(\ka_{i})$ be the eigenvalues of $A$. The eigenvalues of $F'$ can be read off \eqref{DerHk-1}. They are
\eq{F^{i}=\sum_{l=1}^{m}l\fr{\del \psi}{\del p_{l}}\ka_{i}^{l-1}=\fr{\del f}{\del \ka_{i}},}
due to the chain rule.
}

There also follows a representation for the second derivatives of the function $F$. Proofs for the case that $F$ is defined on the subspace of selfadjoint operators with respect to a fixed metric can be found in \cite[Thm.~5.1]{Andrews:/2007}, \cite[Lemma~2.1.14]{Gerhardt:/2006} and \cite{Silhavy:/2000}, where in the latter even higher derivatives are treated.
The proof presented here is by direct differentiation of \eqref{DerHk-2}. It extends similar proofs used in the context of tensor valued functions in \cite{BowenWang:01/1970,BowenWang:01/1971,ChadwickOgden:01/1971,ChadwickOgden:01/1971b} to the case $n>3$ and diagonalisable $A.$
There are several other very recent results \cite{JiangSendov:03/2018}, which address similar questions in the context of operator monotone functions and $k$-isotropic functions. Also compare the comprehensive thesis \cite{Jiang:/2017}, as well as \cite{LewisSendov:/2001} and \cite{Sendov:07/2007}. 

{\prop{\label{2DerHk}
Let $V$ be an $n$-dimensional real vector space and let $F$ and $f$ be as in \cref{DerHk} with $r\geq 2$. Let $A\in \cD_{\G}(V)$ and let $(\eta^{i}_{j})$ be a matrix representation of some $\eta\in \cL(V)$ with respect to a basis of eigenvectors of $A$. Then there holds
\eq{\label{2DerHk-A}d^{2}F(A)(\eta,\eta)=\sum_{i,j=1}^{n}\fr{\del^{2} f}{\del\ka_{i}\del\ka_{j}}\eta^{i}_{i}\eta_{j}^{j}+\sum_{i\neq j}^{n}\fr{\fr{\del f}{\del \ka_{i}}-\fr{\del f}{\del\ka_{j}}}{\ka_{i}-\ka_{j}}\eta^{i}_{j}\eta^{j}_{i},}
where $f$ is evaluated at the $n$-tuple $(\ka_{i})$ of corresponding eigenvalues. The latter quotient is also well defined in case $\ka_{i}=\ka_{j}$ for some $i\neq j$. 
}}

\begin{proof}
Starting from \eqref{DerHk-2} we can calculate for all $A\in \Om\sub\cL(V)$ and $B,C\in \cL(V)$, that
\eq{\label{2DerHk-1} d^{2}F(A)(B,C)&=\sum_{k,l=1}^{m}\fr{\del^{2}\psi}{\del P_{l}\del P_{k}}(d P_{l}(A)B)(d P_{k}(A)C)\\
				&\hp{=}+\sum_{k=1}^{m}\fr{\del \psi}{\del P_{k}}d^{2}P_{k}(A)(B,C).}
From \eqref{DerHk-3} we obtain, already inserting $B=C=\eta=\hat{\eta}+\ti{\eta},$
where $\hat{\eta}$ is the diagonal part of $\eta$ in a basis of eigenvectors for $A$ and $\ti\eta$ is the corresponding off-diagonal part $\ti{\eta}=\eta-\hat{\eta}$,
\eq{\label{2DerHk-3}d^{2}P_{k}(A)(\eta,\eta)&=k\sum_{l=1}^{k-1}\tr(A^{l-1}\circ \eta\circ A^{k-1-l}\circ\eta)\\
					&=k\sum_{l=1}^{k-1}\big(\tr(A^{l-1}\circ \hat\eta\circ A^{k-1-l}\circ\hat\eta)\\
					&\hp{=k\sum\big(}+\tr(A^{l-1}\circ \ti\eta\circ A^{k-1-l}\circ \ti\eta)\big).}
Using the specific basis of eigenvectors we get
\eq{\label{2DerHk-4}d^{2}P_{k}(A)(\eta,\eta)&=k(k-1)\sum_{i=1}^{n}\ka_{i}^{k-2}(\eta^{i}_{i})^{2}+k\sum_{l=1}^{k-1}\sum_{i,j=1}^{n}\ka_{i}^{l-1}\ka_{j}^{k-1-l}\ti\eta^{i}_{j}\ti\eta^{j}_{i}\\
					&=\sum_{i,j=1}^{n}\fr{\del^{2} p_{k}}{\del \ka_{i}\del\ka_{j}}\eta^{i}_{i}\eta^{j}_{j}+\sum_{i\neq j}k\fr{\ka_{i}^{k-1}-\ka_{j}^{k-1}}{\ka_{i}-\ka_{j}}\eta^{i}_{j}\eta_{i}^{j}\\
					&=\sum_{i,j=1}^{n}\fr{\del^{2} p_{k}}{\del \ka_{i}\del\ka_{j}}\eta^{i}_{i}\eta^{j}_{j}+\sum_{i\neq j}\fr{\fr{\del p_{k}}{\del \ka_{i}}-\fr{\del p_{k}}{\del\ka_{j}}}{\ka_{i}-\ka_{j}}\eta^{i}_{j}\eta^{j}_{i}.}
Hence the claimed result holds for the power sums. Returning to \eqref{2DerHk-1} we obtain, also using \cref{EVDerHk},
\eq{d^{2}F(A)(\eta,\eta)&=\sum_{k,l=1}^{m}\fr{\del^{2}\psi}{\del P_{l}\del P_{k}}(d P_{l}(A)\hat\eta)(d P_{k}(A)\hat\eta)+\sum_{k=1}^{m}\fr{\del \psi}{\del P_{k}}d^{2}P_{k}(A)(\eta,\eta)\\
				&=\sum_{i,j=1}^{n}\fr{\del^{2}f}{\del\ka_{i}\del\ka_{j}}\eta^{i}_{i}\eta^{j}_{j}+\sum_{k=1}^{m}\fr{\del \psi}{\del P_{k}}\sum_{i\neq j}\fr{\fr{\del p_{k}}{\del \ka_{i}}-\fr{\del p_{k}}{\del\ka_{j}}}{\ka_{i}-\ka_{j}}\eta^{i}_{j}\eta^{j}_{i},}
from which the claim follows due to the chain rule. Also in this formula, the quotient makes sense even if $\ka_{i}=\ka_{j}$, since the singularity in this fraction is removable, as can be seen from \eqref{2DerHk-4}.
\end{proof}

\begin{rem}
The representation formulae \eqref{EVDerHk-1} and \eqref{2DerHk-A} are only valid a diagonalisable $A$, since their expressions make use of a particular basis of eigenvectors. Formulae which are valid for arbitrary $A\in \Om$ are given, though a little less explicit, in \eqref{DerHk-4} and \eqref{2DerHk-3}. They are still easy enough to serve as a computational tool, particularly in low dimensions.
\end{rem}

Although in the previous proof we have already seen an explicit expression for the quotient term in \eqref{2DerHk-A}, we want to at least mention another representation. It appeared in  \cite[Lemma~2.1.14]{Gerhardt:/2006} and \cite{Silhavy:/2000}, also compare \cite[Lemma~2]{EckerHuisken:02/1989}. The proof is similar to these references.

{\lemma{\label{2DerHkMixed}
Let $f$ be as in \cref{DerHk} with $r\geq 2$ and suppose that $\G$ is convex. Then there holds

\eq{\fr{\fr{\del f}{\del \ka_{i}}-\fr{\del f}{\del \ka_{j}}}{\ka_{i}-\ka_{j}}=\fr 12 \int_{0}^{1}\br{\fr{\del^{2}f}{\del\ka_{i}^{2}}-2\fr{\del^{2}f}{\del\ka_{i}\del\ka_{j}}+\fr{\del^{2}f}{\del\ka_{j}^{2}}},}
 where the integrand is evaluated along the line segment
\eq{\s(t)=\ka+t\fr{\ka_{j}-\ka_{i}}{2}\br{e_{i}-e_{j}}.}
}}

\subsection*{An alternative proof}

Let us have a look at a second nice proof of \cref{2DerHk}, the idea of which appeared in \cite[Lemma~3.2]{Silhavy:/2000}. I owe thanks to the anonymous referee for the observation that this method can also be applied in our situation. It is based on the fact that the function $F$, as given in \cref{DerHk}, is $\mrm{Gl_{n}}(V)$-invariant:
\eq{\label{Inv}F(SAS^{-1})=F(A)\q\fa A\in \cL(V)~ \fa S\in \mrm{Gl_{n}}(V).} 
In \cite[Lemma~3.2]{Silhavy:/2000} this property held for all orthogonal transformations $S$ of a subspace of self-adjoint operators, but the proof basically carries over. Let us repeat it quickly here.

We suppose that all eigenvalues of $A$ are mutually different. The general case can then be treated by approximation as in \cite{Silhavy:/2000}. Differentiating the relation \eqref{Inv} with respect to $A$ in direction of an arbitrary $\eta\in \cL(V)$ we obtain for all $S\in \mrm{Gl_{n}}(V)$, that
\eq{\label{Inv-2}dF(SAS^{-1})(S\eta S^{-1})=dF(A)(\eta).}
In particular, choosing $S=e^{tW}$ for arbitrary $W\in \cL(V)$, $t\in \R$, and differentiating \eqref{Inv-2} with respect to $t$ at $t=0$ gives
\eq{d^{2}F(A)(WA-AW,\eta)=dF(A)(\eta W-W\eta).}
On the other hand, writing 
\eq{\eta=\hat\eta+\ti{\eta},}
with diagonal $\hat\eta$ and off-diagonal $\ti{\eta},$
we have
\eq{d^{2}F(A)(\eta,\eta)=d^{2}F(A)(\hat\eta,\hat\eta)+2d^{2}F(A)(\hat\eta,\ti{\eta})+d^{2}F(A)(\ti{\eta},\ti\eta).}
With respect to a basis of eigenvectors for $A$ and $F'(A)$ we define
\eq{W^{i}_{j}=\fr{\ti{\eta}_{j}^{i}}{\ka_{j}-\ka_{i}},}
which implies
\eq{W^{i}_{k}A^{k}_{j}-A^{i}_{k}W^{k}_{j}=\ti{\eta}^{i}_{j}}
and hence
\eq{d^{2}F(A)(\hat\eta,\ti\eta)=dF(A)(\hat\eta W-W\hat\eta)=0}
and 
\eq{d^{2}F(A)(\ti\eta,\ti\eta)=DF(A)(\ti\eta W-W\ti\eta)=\sum_{i\neq j}\fr{\fr{\del f}{\del\ka_{i}}-\fr{\del f}{\del\ka_{j}}}{\ka_{i}-\ka_{j}}\eta^{i}_{j}\eta^{j}_{i}.}
Finally, since $A$ and $\hat\eta$ are simultaneously diagonal, we have
\eq{d^{2}F(A)(\hat\eta,\hat\eta)&=\fr{d}{dt}\br{dF(A+t\hat\eta)(\hat\eta)}_{|t=0}\\
	&=\fr{d}{dt}\br{\fr{\del f}{\del \ka_{i}}(\ka+t(\eta^{i}_{i}))\eta^{i}_{i}}_{|t=0}\\
	&=\fr{\del^{2}f}{\del\ka_{i}\del\ka_{j}}(\ka)\eta^{i}_{i}\eta_{j}^{j}}
and \cref{2DerHk} follows. \hfill$\square$

There is a slight advantage of the first proof of \cref{2DerHk}, namely that the calculation in \eqref{2DerHk-4} gives a precise description of why the term involving $\ka_{i}-\ka_{j}$ in the denominator also makes sense in case of coalescing eigenvalues.

\section{Functions on bilinear forms}
There is a useful relation of our maps $F\cn \Om\sub\cL(V)\ra \R$ to maps which are defined on bilinear forms. First we need several definitions.

{\defn{Let $V$ be a finite dimensional real vector space.
\begin{itemize}

\item[(i)] We denote the vector space of bilinear forms on $V$ by $\cB(V)$. The space of bilinear forms on the dual space $V^{*}$ is denoted by $\cB^{*}(V)$. The respective subsets of symmetric and positive definite forms will be denoted by $\cB_{+}(V)$ and $\cB^{*}_{+}(V)$.

\item[(ii)] For $a\in\cB(V)$ and $b\in \cB^{*}(V)$ we set
\eq{a_{*}\cn V&\ra V^{*}\\
		v&\mt a(v,\cdot)}
and
\eq{b^{*}\cn V^{*}&\ra V\\
			\phi&\mt J^{-1}\br{b(\phi,\cdot)},}
where $J\cn V\ra V^{**}$ is the canonical identification given by
\eq{v\mt \br{\phi\mt \phi(v)}.}

\item[(iii)] Let $a\in \cB(V)$ and $b\in \cB^{*}(V)$, then we define $b\ast a\in \cL(V)$ by contraction, i.e.
\eq{b\ast a=b^{*}\circ a_{*}.} 
\item[(iv)] For $g\in \cB_{+}(V)$ we define $g^{-1}\in \cB^{*}_{+}(V)$ by requiring
\eq{g^{-1}\ast g=\id.}
\item[(v)] For $a\in \cB(V)$ and $g\in \cB_{+}(V)$ we define the operator $a^{\sharp_{g}}\in \cL(V)$ by
\eq{a^{\sharp_{g}}=g^{-1}\ast a}
\item[(vi)] For any bilinear form $a$ on either $V$ or $V^{*}$ we denote by $\hat{a}$ the symmetrisation, i.e.
\eq{\hat{a}(v,w)=\fr 12\br{a(v,w)+a(w,v)}.} 
\end{itemize}
}}

\begin{rem}
For $a\in \mc{B}(V)$ and $g\in \mc{B}_{+}(V)$ we have
\eq{a(v,w)=g(a^{\sharp_{g}}(v),w)\q\fa v,w\in V.}
\end{rem}

The following construction is very useful.

{\prop{\label{BilF}
Let $V$ be an $n$-dimensional real vector space, $\Om\sub\cL(V)$ open and $F$ be as in \cref{DerHk}. Define
\eq{\Phi\cn \La\sub \cB_{+}(V)\x\cB(V)&\ra \R\\
				(g,h)&\mt F( g^{-1}\ast \hat h), }
where $\La$ is the open subset such that $g^{-1}\ast\hat h\in \Om$ for all $(g,h)\in\La$.
Then $\Phi$ is as smooth as $F$ and the partial derivative of $\Phi$ at $(g,h)$ with respect to $h$ can be regarded as a symmetric bilinear form,
\eq{\fr{\del\Phi}{\del h}(g,h)\in \cB^{*}(V).}
 Furthermore the derivatives of $F$ and $\Phi$ are related by 
\eq{\label{BilF-A}\fr{\del\Phi}{\del h}(g,h)a=\tr(F'(g^{-1}\ast \hat h)\circ\hat a^{\sharp_{g}})=dF(g^{-1}\ast \hat h)\hat a^{\sharp_{g}}.}
}}

\pf{
Since the map $h\mt g^{-1}\ast \hat h$ is linear, we obtain
\eq{\fr{\del\Phi}{\del h}(g,h)a=\tr\br{F'\circ (g^{-1}\ast \hat a) }}
and it can be regarded as a symmetric bilinear form acting on pairs $(\xi,\zeta)$ via letting it act on $\xi\otimes \zeta$.
}

\section{Properties of symmetric functions}

We investigate some special properties associated to symmetric functions, which are particularly related to applications in geometric flows. The most crucial one, the monotonicity, usually ensures that a flow is parabolic. Define
\eq{\G_{+}=\{(\ka_{i})\in \R^{n}\cn \ka_{i}>0\q\fa 1\leq i\leq n\}.}

{\defn{
Let $\G\sub\R^{n}$ open and symmetric, $r\geq 1$ and let $f\in C^{r}(\G)$ be symmetric.
\begin{itemize}
\item[(i)]~$f$ is called {\it{strictly monotone}}, if 
\eq{\fr{\del f}{\del\ka_{i}}(\ka)> 0\q\fa \ka\in\G~\fa 1\leq i\leq n.}
\item[(ii)] Let $\G$ in addition be a cone, then $f$ is called {\it{homogeneous of degree $p\in \R$}} if
\eq{f(\la \ka)=\la^{p}f(\ka)\q\fa \la>0~\fa \ka\in\G.}
\item[(iii)] A nowhere vanishing function $f\in C^{r}(\G_{+})$, $r\geq 2$, is called {\it{inverse concave (inverse convex)}}, if the so-called {\it{inverse symmetric function}} $\ti f\in C^{r}(\G_{+})$, defined by
\eq{\ti f(\ka_{i})=\fr{1}{f(\ka_{i}^{-1})},}
is concave (convex). 
\end{itemize}
}}

These properties carry over to the function $F$ from \cref{DerHk} in the following sense.

{\prop{\label{PropF}
Let $V$ be an $n$-dimensional real vector space, $\G\sub\R^{n}$ open and symmetric, $r\geq 1$ and let $f\in C^{r}(\G)$ and $F\in C^{r}(\Om)$ be as in \cref{DerHk}. 
Then there hold:
\begin{itemize}
\item[(i)] If $f$ is strictly monotone, then $F'(A)$ only has positive eigenvalues at all $A\in \cD_{\G}(V)$ and the bilinear form $\fr{\del\Phi}{\del h}$ from \cref{BilF} is positive definite at all $(g,h)$ with $g^{-1}\ast\hat{h}\in \cD_{\G}(V)$.
\item[(ii)] If $\G$ is a cone and $f$ is homogeneous of degree $p$, then $\cD_{\G}(V)$ is a cone and $F_{|\cD_{\G}(V)}$ is homogeneous of degree $p$.
\item[(iii)]\label{PropF-C} If $r\geq 2$, $\G$ is convex and $f$ is concave, then $F$ satisfies
\eq{d^{2}F(A)(\eta,\eta)\leq 0}
for all $\eta$ having a symmetric matrix representation with respect to a basis of eigenvectors of $A$. The reverse inequality holds if $f$ is convex.
 \end{itemize}
}}

\pf{
(i)~$F'(A)$ has positive eigenvalues due to \cref{EVDerHk}. From \eqref{BilF-A} we obtain (omitting the arguments) for $0\neq \xi\in V$,
\eq{\fr{\del\Phi}{\del h}(\xi,\xi)=\fr{\del\Phi}{\del h}(\xi\otimes\xi)=dF(\xi\otimes\xi)^{\sharp_{g}}>0.}

(ii) Let $A\in \cD_{\G}(V)$ and $\la>0.$ Then the claim follows from $\mrm{EV}(\la A)=\la\mrm{EV}(A)$.

(iii)~Follows immediately from \eqref{2DerHk} and \cref{2DerHkMixed}. 
}

In \cref{PropF-C}, item (iii), the restriction to symmetric $\eta$ is indeed necessary, as can be seen from \cref{example}

The following estimates for $1$-homogeneous resp. inverse concave curvature functions are very useful and are also needed in \cite{BIS4}. The idea for the first statement comes from \cite[Thm.~2.3]{Andrews:/2007} and also appeared in a similar form in \cite[Lemma~14]{BIS1}. The proof for the second statement, however appearing in a slightly different form, can be found in \cite[p.~112]{Urbas:/1991}.

{\prop{\label{InvConc}
Let $V$ be an $n$-dimensional real vector space and $r\geq 1$. Let $f\in C^{r}(\G_{+})$ and $F\in C^{r}(\Om)$ be as in \cref{DerHk} with $f$ symmetric, positive, strictly monotone and homogeneous of degree one. 
Then there hold:
\begin{itemize}
\item[(i)] For every pair $A\in \cD_{\G_{+}}(V)$ and $g\in \cB_{+}(V)$ such that $A$ is self-adjoint with respect to $g$, there holds for all $\eta\in \cL(V)$ that
\eq{dF(A)(\ad_{g}(\eta)\circ A^{-1}\circ\eta)\geq F^{-1}\br{dF(A)\eta}^{2},}
where $\ad_{g}(\eta)$ is the adjoint of $\eta$ with respect to $g$.
\item[(ii)] If $f$ is inverse concave, then for every pair $A\in \cD_{\G_{+}}(V)$ and $g\in \cB_{+}(V)$ such that $A$ is self-adjoint with respect to $g$, there holds
\eq{d^{2}F(A)(\eta,\eta)+2dF(A)(\eta\circ A^{-1}\circ\eta)\geq 2F^{-1}\br{dF(A)\eta}^{2},}
for all $g$-selfadjoint $\eta$.
\end{itemize}
}}

\pf{
(i)~Note that for each $A\in \cD_{\G_{+}}(V)$ the kernel $S$ of the map 
\eq{dF(A)\cn \cL(V)\ra \R}
has dimension $n^2-1$, due to the homogeneity which implies
\eq{dF(A)A=F(A)>0.}
Now let $\eta\in \cL(V)$, then there exists a decomposition
\[\eta=aA+\xi,\]
where $\xi\in S$. Hence, omitting the argument $A$ of $F$,
\eq{dF(\ad_{g}(\eta)\circ A^{-1}\circ\eta)&=adF(\eta)+adF(\ad_{g}(\xi))+dF(\ad_{g}(\xi)\circ A^{-1}\circ\xi)\\
            &\geq adF(\eta),}
since $F'$ and $A$ can be diagonalised simultaneously. The result follows from $F=dF(A)=a^{-1}dF(\eta).$

(ii)~For the inverse symmetric function $\ti f$ the corresponding $\ti F$ has the property
\eq{\label{InvConc-1}\ti F(A)=\fr{1}{F(A^{-1})}\q\fa A\in \cD_{\G_{+}}(V).} 
Thus we may differentiate $\ti F$ using this formula, if we restrict to directions $B$ which are self-adjoint with respect to $g$.  Hence for all $g$-selfadjoint $A\in\cD_{\G_{+}}(V)$ we get
\eq{d\ti F(A)B=\ti F^{2}dF(A^{-1})(A^{-1}\circ B\circ A^{-1})}
and, omitting arguments,
\eq{d^{2}\ti F(B,B)&=2\ti F^{3}\br{dF(A^{-1}\circ B\circ A^{-1})}^{2}\\
			&\hp{=}-\ti F^{2}d^{2}F(A^{-1}\circ B\circ A^{-1},A^{-1}\circ B\circ A^{-1})\\
			&\hp{=}-2\ti F^{2}dF(A^{-1}\circ B\circ A^{-1}\circ B\circ A^{-1}),}
where $\ti F=\ti F(A)$ and $F=F(A^{-1})$. Since $\ti f$ is inverse concave, there holds
\eq{d^{2}\ti F(B,B)\leq 0}
for all $g$-selfadjoint $B$. For some $g$-selfadjoint $\eta$ set
\eq{B=A\circ\eta\circ A} to obtain
\eq{d^{2}F(\eta,\eta)+2dF(\eta\circ A\circ\eta)\geq 2 F^{-1}\br{dF(\eta)}^{2},}
where we again have in mind $F=F(A^{-1}).$ The result follows.
}

\section{Examples}
Let us have a look at some familiar symmetric functions, their corresponding associated operator functions and their properties. The most important examples are the elementary symmetric polynomials satisfying
\eq{s_{k}\circ\mrm{EV}(A)=\fr{1}{k!}\fr{d^{k}}{dt^{k}}\det(I+tA)_{|t=0},}
compare \cite[equ.~(2.1.31)]{Gerhardt:/2006}. $s_{k}$ is strictly monotone on the set
\eq{\G_{k}=\{\ka\in\R^{n}\cn s_{1}(\ka)>0,\dots,s_{k}(\ka)>0\},}
which is equal to the connected component of the set $\{s_{k}>0\}$ containing $\G_{+},$ compare \cite[Prop.~2.6]{HuiskenSinestrari:09/1999}. Obviously $s_{1}$ is also concave and convex.

Define the quotients
\eq{q_{k}\cn \G_{k-1}&\ra \R\\
				q_{k}&=\fr{s_{k}}{s_{k-1}}.}
These are homogeneous of degree one and concave, cf. \cite[Thm.~2.5]{HuiskenSinestrari:09/1999}. On $\G_{+}$ the $q_{k}$ are also strictly monotone and inverse concave, cf. \cite[Thm.~2.6]{Andrews:/2007}.
Also the functions
\eq{f=\br{\fr{s_{k}}{s_{l}}}^{\fr{1}{k-l}}, \q 0\leq l<k\leq n,}
share all these properties on $\G_{+}$, \cite[p.~23]{Andrews:/2007}. More examples of such curvature functions can be found in \cite{Andrews:/2007}.

\section{Loss of regularity}
In this final section we discuss the regularity properties of the associated operator function $F$ and show be means of an example that the loss of regularity from $f$ to $\psi$ in the correspondence
\eq{f=\psi(p_{1},\dots,p_{m})}
also leads, in general, to the same loss of regularity from $f$ to 
\eq{F\cn \Om\ra \R} in the relation
\eq{\label{Loss-1}F_{|\cD_{\G}(V)}=f\circ\mrm{EV}_{|\cD_{\G}(V)}.}
Consider the following example:
\eq{f(\ka_{1},\ka_{2})=(\ka_{1}^{2}+\ka_{2}^{2})^{\fr 32}.}
Then $f\in C^{2}(\R^{2})$.
Since $F$ is required to satisfy \eqref{Loss-1} and the {\it{open}} domain $\Om$ of $F$ has to contain the zero matrix,
we must use 
\eq{\psi(x_{1},\dots,x_{m})=|x_{2}|^{\fr 32}}
to connect to $f$ (note that $P_{2}(A)$ can be negative). Hence
\eq{F\cn \cL(V)&\ra \R\\ 
	F(A)&=\psi(P_{2}(A))=|\tr(A^{2})|^{\fr 32}}
is an associated operator function.
Writing, with respect to a basis,
\eq{A=\begin{pmatrix} w & x\\
				y & z \end{pmatrix},}
we see that
\eq{F(A)=F(w,x,y,z)=|w^{2}+2xy+z^{2}|^{\fr 32},}
which is not $C^{2}$, since its restriction to the straight line
\eq{\label{Loss-2}x\mt (0,x,1,0)}
is not $C^{2}$. It is in fact only as smooth as $\psi$. This is in sharp contrast to the regularity of the restriction to a subspace of $g$-selfadjoint operators,
\eq{F\cn \Sigma_{g}(V)\ra \R,}
which has the same regularity as $f$, cf. \cite{Ball:/1984,Silhavy:/2000}. The crucial difference is that the variations in \eqref{Loss-2} are not allowed, since one must remain within the class of symmetric matrices.

\section*{Acknowledgements}
I would like to thank the anonymous referee for his/her interest in this work, for pointing out the second proof of \cref{2DerHk} and for the very helpful comments which made the exposition much more complete.
I would also like to thank Prof. Guofang Wang for a helpful discussion on isotropic functions.

\bibliographystyle{hamsplain}
\bibliography{bibliography}

\providecommand{\bysame}{\leavevmode\hbox to3em{\hrulefill}\thinspace}
\providecommand{\href}[2]{#2}
\begin{thebibliography}{10}

\bibitem{Andrews:/2007}
Ben Andrews, \emph{Pinching estimates and motion of hypersurfaces by curvature
  functions}, J. Reine Angew. Math. \textbf{608} (2007), 17--33.

\bibitem{Ball:/1984}
John Ball, \emph{Differentiability properties of symmetric and isotropic
  functions}, Duke Math. J. \textbf{51} (1984), no.~3, 699--728.

\bibitem{Barbancon:/1972}
G\'erard Barban\c{c}on, \emph{Th\'eor\`eme de {N}ewton pour les fonctions de
  classe ${C}^r$}, Ann. Sci. {\'E}c. Norm. Sup{\'e}r. (4) \textbf{5} (1972),
  no.~3, 435--457.

\bibitem{BowenWang:01/1970}
Ray Bowen and Chao~Chen Wang, \emph{Acceleration waves in inhomogeneous
  isotropic elastic bodies}, Arch. Rat. Mech. Anal. \textbf{38} (1970), no.~1,
  13--45.

\bibitem{BowenWang:01/1971}
\bysame, \emph{Corrigendum: Acceleration waves in inhomogeneous isotropic
  elastic bodies}, Arch. Rat. Mech. Anal. \textbf{40} (1971), no.~5, 403--403.

\bibitem{BIS1}
Paul Bryan, Mohammad~N. Ivaki, and Julian Scheuer, \emph{Harnack inequalities
  for evolving hypersurfaces on the sphere}, To appear in Commun. Anal. Geom.,
  {\href{https://arxiv.org/abs/1512.03374}{arxiv:1512.03374}}, 2015.

\bibitem{BIS4}
\bysame, \emph{Harnack inequalities for curvature flows in {R}iemannian and
  {L}orentzian manifolds}, preprint,
  {\href{https://arxiv.org/abs/1703.07493}{arxiv:1703.07493}}, 2017.

\bibitem{ChadwickOgden:01/1971}
Peter Chadwick and Raymond Ogden, \emph{On the definition of elastic moduli},
  Arch. Rat. Mech. Anal. \textbf{44} (1971), no.~1, 41--53.

\bibitem{ChadwickOgden:01/1971b}
\bysame, \emph{A theorem of tensor calculus and its application to isotropic
  elasticity}, Arch. Rat. Mech. Anal. \textbf{44} (1971), no.~1, 54--68.

\bibitem{EckerHuisken:02/1989}
Klaus Ecker and Gerhard Huisken, \emph{Immersed hypersurfaces with constant
  {W}eingarten curvature}, Math. Ann. \textbf{283} (1989), no.~2, 329--332.

\bibitem{Gerhardt:/2006}
Claus Gerhardt, \emph{Curvature problems}, Series in Geometry and Topology,
  vol.~39, International Press of Boston Inc., Sommerville, 2006.

\bibitem{Glaeser:01/1963}
Georges Glaeser, \emph{Fonctions compos\'ees diff\'erentiables}, Ann. Math.
  \textbf{77} (1963), no.~1, 193--209.

\bibitem{HuiskenSinestrari:09/1999}
Gerhard Huisken and Carlo Sinestrari, \emph{Convexity estimates for mean
  curvature flow and singularities of mean convex surfaces}, Acta Math.
  \textbf{183} (1999), no.~1, 45--70.

\bibitem{Jiang:/2017}
Tianpei Jiang, \emph{Properties of k-isotropic functions}, Ph.D. thesis,
  University of Western {O}ntario, 2017.

\bibitem{JiangSendov:03/2018}
Tianpei Jiang and Hristo Sendov, \emph{A unified approach to operator monotone
  functions}, Linear Algebra Appl. \textbf{541} (2018), 185--210.

\bibitem{LewisSendov:/2001}
Adrian Lewis and Hristo Sendov, \emph{Twice differentiable spectral functions},
  SIAM J. Matrix. Anal. Appl. \textbf{23} (2001), no.~2, 368--386.

\bibitem{Mead:10/1992}
David Mead, \emph{Newton's identities}, Am. Math. Mon. \textbf{99} (1992),
  no.~8, 749--751.

\bibitem{Sendov:07/2007}
Hristo Sendov, \emph{The higher-order derivatives of spectral functions},
  Linear Algebra Appl. \textbf{424} (2007), no.~1, 240--281.

\bibitem{Silhavy:/2000}
Miroslav Silhav\'y, \emph{Differentiability properties of isotropic functions},
  Duke Math. J. \textbf{104} (2000), no.~3, 367--373.

\bibitem{Urbas:/1991}
John Urbas, \emph{An expansion of convex hypersurfaces}, J. Differ. Geom.
  \textbf{33} (1991), no.~1, 91--125.

\bibitem{Zedek:/1965}
Mishael Zedek, \emph{Continuity and location of zeros of linear combinations of
  polynomials}, Proc. Am. Math. Soc. \textbf{16} (1965), no.~1, 78--84.

\end{thebibliography}

\end{document}